\RequirePackage{fix-cm}
\documentclass[smallextended]{svjour3}       
\smartqed  

\usepackage{graphicx}
\usepackage{amssymb,amsmath,amscd,amssymb,amsxtra}
\usepackage{epsfig}
\usepackage{mathptmx}
\usepackage{helvet}
\usepackage{courier}
\usepackage{txfonts}
\usepackage{tikz} 
\usetikzlibrary{arrows}
\usepackage{type1cm} 
\usepackage{esvect}
\usepackage{enumitem}
\usepackage{titlesec}
\usepackage{float}
\usepackage{bm}
\usepackage{comment}
\usepackage{xcolor}
\titlelabel{\thetitle.\quad}

\newcommand{\diam} {\mathop{\rm diam}\nolimits}
\newcommand{\dr} {\mathop{\rm degrig}\nolimits}

\newcommand{\Cl}  {\mathop{\rm Cl}\nolimits}

\newcommand{\Mod}  {\mathop{\rm mod}\nolimits}
\newcommand{\degrig}  {\mathop{\rm degrig}\nolimits}
\newcommand{\Hed}  {\mathop{\rm Hed}\nolimits}
\newcommand{\HR}  {\mathop{\rm HR}\nolimits}

\newtheorem{observation}[theorem]{Observation}
\newtheorem{construction}[theorem]{Construction}
\spnewtheorem{notation}[theorem]{Notation}{\bfseries}{\itshape}

\newcommand{\NN}{\mathbb N}

\usetikzlibrary {decorations.pathmorphing}
\tikzset{snake it/.style={decorate, decoration=snake}}
\tikzset{step it/.style={decorate, decoration={random steps,segment length=3pt,amplitude=1.5pt}}}
\tikzset{star it/.style={decorate, decoration={random steps,segment length=2pt,amplitude=0.5pt}}}

\begin{document}
\def\joinrel{\mkern-3mu}
\newcommand{\varproj}{\displaystyle \lim_{\multimapinv\joinrel-\joinrel-}}

\title{Rigid continua containing an arbitrary Peano continuum  \thanks{This work is supported in part by the Slovenian Research and Innovation Agency (research programs P1-0403, J1-70016 and SN-ZRD/22-27/0552).}
}


\author{Matev\v z \v Crepnjak \and				
				Teja Kac \and Tina Sovi\v c
}


\institute{M. \v Crepnjak \at
			Faculty of Chemistry and Chemical Engineering, University of Maribor, Smetanova 17, SI-2000 Maribor, Slovenia \\
              Faculty of Natural Sciences and Mathematics, University of Maribor, Koro\v{s}ka 160, SI-2000 Maribor,
   Slovenia \\
              \email{matevz.crepnjak@um.si}           
           \and
           T. Kac \at
              Faculty of Natural Sciences and Mathematics, University of Maribor, Koro\v{s}ka 160, SI-2000 Maribor,
   Slovenia \\
              \email{teja.kac1@um.si}    
             \and
           T. Sovi\v c \at
              Faculty of Civil Engineering, Transportation Engineering and Architecture, University of Maribor, Smetanova 17, SI-2000 Maribor, Slovenia \\
              Faculty of Natural Sciences and Mathematics, University of Maribor, Koro\v{s}ka 160, SI-2000 Maribor,
   Slovenia \\
              \email{tina.sovic@um.si} 
}

\date{Received: date / Accepted: date}

\maketitle

\begin{abstract}

Bani\v c and Kac, and later \v Crepnjak and Kac, constructed an uncountable family of continua with different types of rigidity. For every continuum $X$ with the mentioned property, the inverse limit $\varprojlim\{X,f_n\}$ is homeomorphic to $X$ for an arbitrary sequence of continuous surjective functions $f_n: X \to X$. 

We extend these results by proving that for any Peano continuum $X$ there exists an uncountable family ${\cal F}$ of rigid, for each positive integer $m$, $\frac{1}{m}$-rigid, and $0$-rigid continua such that each $Y \in {\cal F}$ contains a homeomorphic copy of $X$ and for each $Y \in {\cal F}$  it holds that the inverse limit $\varprojlim\{Y,f_n\}$ is homeomorphic to $Y$ for arbitrarily sequence of continuous surjective functions $f_n: Y \to Y$.

\end{abstract}


\section{Introduction}

In 1955, de Groot posed the question of whether there exists a connected metric space for which no surjective continuous function maps the space onto its proper non-degenerate subspace, \cite{Groot}. Later, Anderson in \cite{Anderson} asked the following: Does there exist a non-degenerate continuum $X$ where the only continuous self-maps $X \to X$ are the identity and constant functions? Moreover, if such a continuum exists, do all of its non-degenerate subcontinua share this property? In  \cite{AndersonChoquet}, Anderson and Choquet gave an example of a planar continuum where each two non-degenerate different subcontinua are non-homeomorphic. In \cite{cook}, Cook constructed two examples of continua $X$ and $Y$ with the following properties:
\begin{itemize}
    \item The continuum $X$ has the property that the identity is the only continuous non-constant self-map $X \to X$.
    \item The continuum $Y$ has the property that for any pair of non-degenerate different subcontinua, there exists no continuous non-constant function from one to the other. 
\end{itemize}
Later, many authors named continua with the second property from above Cook continua.  In this paper, we adopt the same definition. 

In the papers, various types of rigid continua have been introduced (e.g.\ \cite{banicKac,CharalambousKrzempek,crepnjakKac,veronika,jan,Renska,Terasawa,Terasawa2}).
Recently, new types of rigidity of continua have been introduced and, by this, the theory of rigid continua has been extended. In the present paper, we focus on the type introduced in \cite{banicKac,crepnjakKac} where the authors defined $\frac{1}{m}$-rigid and $0$-rigid continua and constructed examples of them. They also partially answered the following question: Find all non-degenerate continua $X$ such that $\varprojlim \{X,f_n\}$ is homeomorphic to $X$ for any sequence $(f_n)$ of continuous surjections $f_n : X \to X$. The first continuum discovered with this property was the pseudoarc. Given the current examples, a complete classification of such continua appears not to be viable. Therefore, we approach this problem from a different perspective. More precisely, for each non-degenerate Peano continuum $X$ we construct a homeomorphic copy $h(X)$ of $X$ and a rigid, $\frac{1}{m}$-rigid, and $0$-rigid continuum $Y$, respectively, such that $h(X) \subseteq Y$ and for each sequence of continuous surjective functions $(f_n)$, $f_n: Y \to Y$ for each positive integer $n$, the following hold:
\begin{itemize}
    \item $f_n|_{h(X)}=1_{h(X)}$ for each positive integer $n$,
    \item $X$ and $h(X)$ are homeomorphic to $\varprojlim\{h(X),f_n|_{h(X)}\}$,
    \item $Y$ is homeomorphic to $\varprojlim\{Y,f_n\}$.
\end{itemize}

We proceed as follows. In Section \ref{sec:preliminaries}, the definitions and notation
 needed in the paper are given. In Section \ref{sec:construction_hed}, we give a general construction of a so-called Hedgehog of $X$, or Hed$(X)$ for short, which is a continuum containing a homeomorphic copy of $X$ and possesses the required properties. In Section \ref{sec:rig->hed}, using some modification of a general construction, we obtain a rigid, $0$-rigid and, for each positive integer $m \geq 2$, also a $\frac{1}{m}$-rigid continuum with desired properties. In Section \ref{sec:application} we give a short note about the application of our results to inverse limits.

\section{Preliminaries}\label{sec:preliminaries}

\begin{definition}
    A {\em continuum} is a non-empty, compact, and connected metric space.
\end{definition}

\begin{notation}
Let $n$ be a positive integer,  $X$ a continuum, and $f:X\to X$ a function. We use $f^n$ to denote the composition: 
$$
f^n=\underbrace{f\circ f\circ\cdots\circ f}_n. 
$$
\end{notation}

\begin{notation}
    With $1_X$ we denote the identity function on a continuum $X$; a function $1_X:X\to X$ is defined by $1_X(x)=x$ for each $x\in X$.
\end{notation}

\begin{definition}
    A {\em Peano continuum} is a locally connected continuum.
\end{definition}

We use the following well-known characterization of Peano continua.

\begin{theorem}\cite[Theorem 8.18]{nadler}\label{Hahn}
    A continuum $X$ is a Peano continuum if and only if $X$ is a continuous image of $[0,1]$.
\end{theorem}

\begin{definition}
    Let $X$ be any continuum and $x_0\in X$. The point $x_0$ is a {\em cut point for $X$} if $X\setminus \{x_0\}$ is not connected. Otherwise, $x_0$ is a {\em non-cut point for $X$}.
\end{definition}

In the next definition, we present Cook continua, whose existence was proved by Cook in \cite{cook}.

\begin{definition}\label{def:Cook_continuum}
    A non-degenerate continuum $X$ is a {\em Cook continuum} if for any non-degenerate subcontinuum $S$ of $X$ and any non-constant continuous function $f: S\to X$, it follows that $$f(x)=x$$ for each $x\in S$.
\end{definition}

The next result can be found in \cite{banicKac}.

\begin{lemma}\cite[Lemma 3.2]{banicKac}\label{lemma:cut_point_Cook}
    Let $X$ be a Cook continuum, and let $x_0\in X$ be any point. Then $x_0$ is a non-cut point for $X$.
\end{lemma}

The next lemma is proved in \cite{snoha}.

\begin{lemma}\cite[Lemma 2.12]{snoha}\label{lemma:cook_funkcije}
    \begin{enumerate}
        \item Let $X$ be a Cook continuum and $S$ be a non-degenerate subcontinuum of $X$. Then $S$ is also a Cook continuum.
        \item If two Cook continua $X_1$ and $X_2$ are homeomorphic, then they are uniquely homeomorphic, i.e., there exists a unique homeomorphism $\varphi: X_1 \to X_2$. Moreover, if $\psi :  X_1 \to X_2$ is a continuous surjection,
then $\psi = \varphi$.
    \end{enumerate}
\end{lemma}

The next result can be found in \cite{crepnjakKac}.
\begin{proposition}\cite[Proposition 1]{crepnjakKac}\label{prop1}
    Let $X$ be any Cook continuum and let $S$ be any non-degenerate subcontinuum of $X$, such that $S \neq X$. If $f: S \to X$ is a continuous function, then $f$ is a constant function.
\end{proposition}

\begin{definition}\label{def:rigid}
    Let $n$ be a positive integer. A continuum $X$ is {\em $\frac{1}{n}$-rigid}, if $$f^n=1_X$$ for any continuous surjective function $f:X\to X$. If a continuum $X$ is not $\frac{1}{n}$-rigid for any positive integer $n$, we say that $X$ is {\em $0$-rigid}. If a continuum $X$ is $\frac{1}{1}$-rigid, then we also say that $X$ is {\em rigid}.
\end{definition}

Note that for each positive integer $n$ and for any $\frac{1}{n}$-rigid continuum $X$ it holds that $X$ is also $\frac{1}{kn}$-rigid for each positive integer $k$. Thus, we define the degree of rigidity of a continuum. 

\begin{definition}
    Let $X$ be a continuum. We define {\em the degree of rigidity}, $\degrig(X)$, of continuum $X$ as follows
    $$\degrig(X)=\min\left(\left\{n\in\NN\ |\ X \text{ is } \frac{1}{n} \text{-rigid}\right\} \cup \{\infty\}\right).$$
\end{definition}

\begin{observation}
    Any Cook continuum has the degree of rigidity $1$ and it is rigid.
\end{observation}

For each positive integer $n$, the existence of a continuum whose degree of rigidity is equal to $n$ is proved in \cite{banicKac}. One can easily see that an arc has the degree of rigidity equal to $\infty$ (i.e.\ an arc is $0$-rigid).

\begin{definition}
An inverse sequence is a double sequence $\{X_n,f_n\}_{n=1}^{\infty}$ of metric spaces $X_n$ and functions $f_n:X_{n+1}\rightarrow X_n$.
The spaces $X_n$ are called the coordinate spaces, and the functions $f_n$ are called the bonding functions.
The inverse limit of an inverse sequence $\{X_n,f_n\}_{n=1}^{\infty}$ is the subspace of the product space $\prod_{n=1}^{\infty}X_n$ that consists of all points $(x_1,x_2,x_3,\ldots)\in \prod_{n=1}^{\infty}X_n$ such that $x_{n+1}=f_n(x_n)$ for each positive integer $n$:
$$
\varprojlim\{X_n,f_n\}_{n=1}^{\infty}=\left\{(x_1,x_2,x_3,\ldots)\in \prod_{n=1}^{\infty}X_n \ | \ x_{n}=f_n(x_{n+1}) \textup{ for each } n\right\}.
$$
\end{definition}

\begin{definition}
Let $n$ be a positive integer and let $(x_1,x_2,\ldots,x_n,\ldots)\in \displaystyle \prod_{n=1}^\infty [0,1]$. A function $\pi_n:\displaystyle \prod_{n=1}^\infty [0,1]\to [0,1]$ is called the {\em $n$-th projection} and it is defined as $$\pi_n(x_1,x_2,\ldots,x_n,\ldots)=x_n.$$
\end{definition}

\begin{definition}
A metric $d$ on a Hilbert cube $\displaystyle \prod_{n=1}^{\infty} [0,1]$ is defined as
$$d(x, y) = \sum_{n=1}^{\infty} \frac{|x_n-y_n|}{2^n},$$
for each $x = (x_1,x_2,x_3,\ldots), y = (y_1,y_2,y_3,\ldots) \in \displaystyle \prod_{n=1}^{\infty} [0,1]$.
\end{definition}

\section{Construction of a fan of Cook continua}\label{sec:construction_fan}

Later in the paper we use so-called rigid stars of Cook continua. In this section, we define these special fans of Cook continua and provide some of their properties.

\begin{construction}\label{construction:rigid_star}
    Let $(A_n)_{n=1}^\infty$ be a sequence of subsets of positive integers such that
\begin{enumerate}
\item $\displaystyle \bigcup_{n=1}^\infty A_n$ is equal to the set of all positive integers,
\item for all positive integers $m$ and $n$ it holds that: $m\neq n$ if and only if $A_m\cap A_n=\emptyset$, and
\item $|A_n|=\infty$ for each positive integer $n$.
\end{enumerate}
Now, for each positive integer $n$, let $Q_n\subseteq \displaystyle \prod_{n=1}^\infty [0,1]$ such that
$$\pi_k(Q_n)=\begin{cases}
\{0\}& \text{;} \quad k\notin A_n,\\
[0,\frac{1}{n}] & \text{;} \quad k\in A_n,\\
\end{cases}$$
for each positive integer $k$.

Let $C$ be a Cook continuum and let $(X_n)$ be a sequence of non-degenerate continua, such that
\begin{itemize}
\item $X_1=C$,
\item $X_{n+1}$ is a subcontinuum of $X_n$ for each positive integer $n$, and
\item $\lim\limits_{n\to\infty}\diam (X_n)=0$.
\end{itemize}
By Lemma \ref{lemma:cook_funkcije}, $X_n$ is a Cook continuum for each positive integer $n$. The existence of such a sequence of continua follows from the Boundary Bumping Theorem, \cite[Theorem 5.4]{nadler}. 

Now, let $\{t_0\}=\displaystyle \bigcap_{n=1}^\infty X_n$.

For each positive integer $n$, let $H_n:X_n\to Q_n$ be an embedding such that $H_n(t_0)=(0,0,\ldots)$. Such embeddings do exist since $ \displaystyle \prod_{n=1}^\infty [0,1]$ is a homogeneous continuum. 

Let $X=\displaystyle \bigcup_{n=1}^\infty H_n(X_n)$.
Then $X$ is called a {\em fan of Cook continua}.
\end{construction}

In this section we denote $\overline{0} = (0,0,\ldots)$.

\begin{observation}\label{obs:X_cont_cut_point}
Any fan of Cook continua is a continuum. Moreover, a point $\overline{0}$ is the only cut point of such a continuum since Cook continua do not have cut points, see Lemma \ref{lemma:cut_point_Cook}.
\end{observation}

\begin{example}
    Any simple fan of Cook continua defined in \cite[Construction 1]{crepnjakKac} is a fan of Cook continua.
\end{example}

The following lemma is a similar result to Lemma 4 in \cite{crepnjakKac}.

\begin{lemma}\label{lemma:rigid_star_0}
    Let $X$ be a fan of Cook continua, and $f:X\to X$ be a continuous surjection. Then $$f(\overline{0})=\overline{0}.$$
\end{lemma}
\begin{proof}
    Suppose that $f(\overline{0})\neq\overline{0}$. Then there exists a positive integer $n$ such that $f(\overline{0}) \in H_n(X_n)$. Let
    $g: X \to H_n(X_n)$ such that
    $$g(x) = \begin{cases}
        x \quad \text{;} & x \in H_n(X_n), \\
        \overline{0} \quad \text{;} & x \notin H_n(X_n),
    \end{cases}$$
    for each $x \in X$. Note that $g$ is a continuous surjective function. Therefore $H_n^{-1} \circ g \circ f \circ H_i:X_i \to X_n$ is a continuous function for each positive integer $i$. Since $X$ is a Cook continuum with $X_i,X_n\subseteq X$ and $(H_n^{-1} \circ g \circ f \circ H_i)(t_0) \neq t_0$, it follows that $H_n^{-1} \circ g \circ f \circ H_i$ is a constant function, which maps $X_i$ to $H_n^{-1}(f(\overline{0}))$.

    If there exists $x\in H_i(X_i)$ such that $f(x)\neq f(\overline{0})$, then $(g\circ f)(\overline{0}) \neq f(\overline{0})$. Therefore, $(H_n^{-1} \circ g \circ f \circ H_i)(H_i^{-1}(x)) \neq H_n^{-1}(f(\overline{0}))$--a contradiction.

    Consequently, $f|_{H_i(X_i)}$ is a constant function for each positive integer $i$--a contradiction, since $f$ is a surjective function.
    \qed
\end{proof}

In the next definition we define a special type of fans of Cook continua.

\begin{definition}
    Let $X$ be a fan of Cook continua such that $X_{n+1}$ is a proper subcontinuum of $X_n$ for each positive integer $n$.
    Then we call $X$ {\em a rigid star of Cook continua}.
\end{definition}

\begin{observation}\label{obs:simple_fan}
    Let $(X_n)$ be a sequence of Cook continua, and let $X$ be a rigid star of Cook continua obtained from $(X_n)$. In \cite{crepnjakKac}, a construction of a simple fan of Cook continua is given in such a way that each $X_n$ is used exactly $n$ times. More precisely, we construct $n$ copies of $X_n$, and we denote them by $X_n^k$ for each positive integer $n$ and for each $k\in\{1,2,\ldots,n\}$ (for all the details see \cite[Construction 1]{crepnjakKac}). 
    
    Now we construct a simple fan of Cook continua by $\bigcup_{n=1}^\infty\bigcup_{k=1}^n X_n^k$. We denote such continuum by $Y$. Continuum $X$ can be embedded in $Y$ in a way that the image of each $H_n(X_n)$ is equal to the $X_n^1$ for each positive integer $n$. So, we understand $X$ as a subcontinuum of $Y$ or $Y$ as an extension of $X$. We use this fact in the next theorem.
\end{observation}

\begin{theorem}\label{thm:rigid_star_is_rigid}
    Any rigid star of Cook continua is rigid.
\end{theorem}

\begin{proof}
    Let $X$ be a rigid star of Cook continua and let $f:X\to X$ be a continuous surjection. From Lemma \ref{lemma:rigid_star_0} it follows that $f(\overline{0})=\overline{0}$. Let $Y$ be a simple fan of Cook continua obtained by a sequence of Cook continua $(X_n)$ as in Observation \ref{obs:simple_fan}. Without loss of generality, we treat the continuum $X$ as a subcontinuum of $Y$.

    We define a function $g:Y\to Y$ in the following way
    $$g(x)=\begin{cases}
    f(x)&  \text{;} \quad x\in X,\\
    x &  \text{;}\quad x\in Y\setminus X,
    \end{cases}$$
    for each $x\in Y$. The function $g$ is a continuous surjection. From \cite[Theorem 1]{crepnjakKac} follows that $g(x)=x$ for each $x\in Y$. Therefore, $f(x)=x$ for each $x\in X$ and thus, $X$ is rigid.
    \qed
\end{proof}

\section{Construction of the hedgehog $\text{Hed}(X)$ of a given continuum $X$}\label{sec:construction_hed}

In this section, we give the first main result of the paper. More precisely, for any Peano continuum $X$ we construct an uncountable family ${\cal F}$ of pairwise non-homeomorphic continua such that for each $Y \in {\cal F}$ it holds that $Y$ is rigid and it contains a homeomorphic copy of $X$. 


\begin{construction}\label{construction:hedgehog}
Let $X$ be an arbitrary continuum. Without loss of generality, suppose $\diam(X) \leq 1$. Let $P=\{p_n\in X \ | \ n \in \mathbb{N} \}$ be a countable dense subset of $X$, for which it holds $p_i\neq p_j$ if and only if $i\neq j$. Now define $h : X \to Q$, where $Q$ is the Hilbert cube, by 
$$h(x) = (d(x,p_1),0,d(x,p_2),0,d(x,p_3),0,\ldots)$$
for each $x \in X$. It follows that $h$ is an injective continuous function. Since $X$ is a continuum and $Q$ is a metric space, $h$ is closed. Consequently, $h$ is an embedding.

Let $(A_i)$ be a sequence of sets such that
\begin{enumerate}
    \item for each positive integer $i$, $A_i \subseteq \{2n \ | \ n \in \mathbb{N} \}$ and $|A_i| = \infty$,
    \item for each positive integers $i$ and $j$ it holds that $A_i \cap A_j = \emptyset$ if and only if $i \neq j$,
    \item $\bigcup_{i=1}^{\infty} A_i = \{2n \ | \ n \in \mathbb{N}\}$.
\end{enumerate}

Let $R = h(P)$ and denote $R = \{\bm{r}_n \ |~ \bm{r}_n=h(p_n) \textup{ for each} \ n \in \mathbb{N}\}$. It follows $R$ is a dense subset of $h(X)$ satisfying the condition $\bm{r}_n \neq \bm{r}_m$ if and only if $n\neq m$. Let $\bm{r}_n = (r_1^n,r_2^n,r_3^n,\ldots)$ for each positive integer $n$. Note that $r_{2i}^n = 0$ for each positive integers $i$ and $n$.

For each positive integer $n$ we define $Q_n\subseteq Q$ such that the following holds
\begin{enumerate}
    \item $\pi_{2k-1}(Q_n) = \{r_{2k-1}^n\}$ for each positive integer $k$,
    \item $\pi_k(Q_n) = [0,1]$ for each $k \in A_n$,
    \item $\pi_{k}(Q_n) = \{0\}$ for each $\displaystyle k \in \bigcup_{j=1, j \neq n} ^{\infty} A_j$. 
\end{enumerate}
It follows that each $Q_n$ is homeomorphic to the Hilbert cube. Note that $Q_n \cap h(X) = \{\bm{r}_n\}$ for each positive integer $n$.

Let $(X_n)$ be a sequence of continua and $(x_n)$ be a sequence of points such that $x_n\in X_n$ for each positive integer $n$.
Now we embed $X_n$ into $Q_n$ by $h_n: X_n \to Q_n$ with the following properties
\begin{enumerate}
    \item $h_n(x_n)=\bm{r}_n$,
    \item $\diam (h_n(X_n)) \leq \frac{1}{n}$.        
\end{enumerate}
Such an embedding exists since the Hilbert cube is a homogeneous continuum. 
Finally, we define the {\em hedgehog of $X, (X_n),(x_n)$} by
$$\Hed(X, (X_n),(x_n)) = h(X) \cup \bigcup_{n=1}^{\infty} h_n(X_n).$$
For $\Hed(X, (X_n),(x_n))$, we also use the simplified notation $\Hed(X)$ and refer to it simply as the {\em hedgehog of $X$}.
\end{construction}







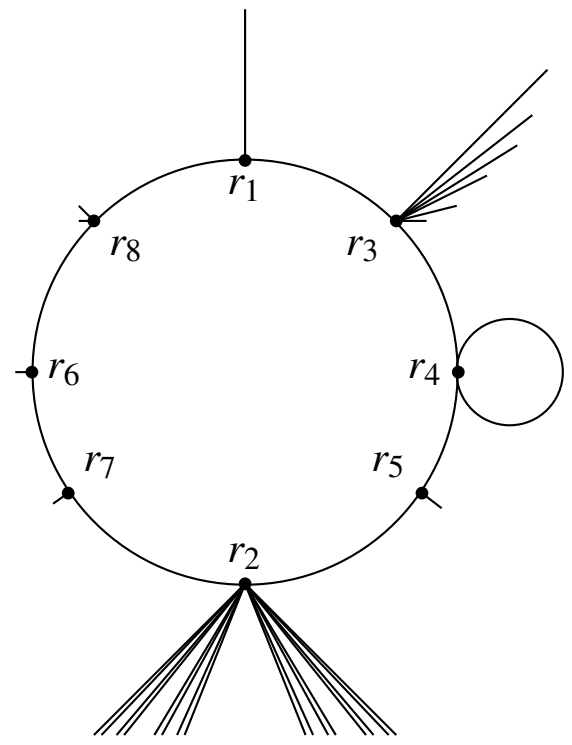
\begin{figure}[H]
    \centering
    \begin{tikzpicture}
        \draw[black] (0,0) circle (40pt);

        \draw[black, thin] (0,1.4) -- (0,2.4);
        \filldraw[black] (0,1.4) circle (1pt) node[anchor=north]{$r_1$};

        \draw[black, thin] (0,-1.4) -- (-1,-2.4);
        \draw[black, thin] (0,-1.4) -- (1,-2.4);
        \draw[black, thin] (0,-1.4) -- (-0.4,-2.4);
        \draw[black, thin] (0,-1.4) -- (0.4,-2.4);
        \draw[black, thin] (0,-1.4) -- (-0.6,-2.4);
        \draw[black, thin] (0,-1.4) -- (0.6,-2.4);
        \draw[black, thin] (0,-1.4) -- (-0.8,-2.4);
        \draw[black, thin] (0,-1.4) -- (0.8,-2.4);
        \draw[black, thin] (0,-1.4) -- (-0.95,-2.4);
        \draw[black, thin] (0,-1.4) -- (0.95,-2.4);
        \draw[black, thin] (0,-1.4) -- (-0.85,-2.4);
        \draw[black, thin] (0,-1.4) -- (0.85,-2.4);
        \draw[black, thin] (0,-1.4) -- (-0.45,-2.4);
        \draw[black, thin] (0,-1.4) -- (0.45,-2.4);
        \draw[black, thin] (0,-1.4) -- (-0.55,-2.4);
        \draw[black, thin] (0,-1.4) -- (0.55,-2.4);
        \filldraw[black] (0,-1.4) circle (1pt) node[anchor=south]{$r_2$};

        \draw[black, thin] (1,1) -- (2,2);
        \draw[black, thin] (1,1) -- (1.9,1.7);
        \draw[black, thin] (1,1) -- (1.8,1.5);
        \draw[black, thin] (1,1) -- (1.6,1.3);
        \draw[black, thin] (1,1) -- (1.4,1.1);  
        \draw[black, thin] (1,1) -- (1.2,1);
        \filldraw[black] (1,1) circle (1pt) node[anchor=north east]{$r_3$};

        \draw[black] (1.75,0) circle (10pt);
        \filldraw[black] (1.41,0) circle (1pt) node[anchor=east]{$r_4$};

        \draw[black, thin] (1.17,-0.8) -- (1.3,-0.9);
        \filldraw[black] (1.17,-0.8) circle (1pt) node[anchor=south east]{$r_5$};
        
        \draw[black, thin] (-1.41,0) -- (-1.52,0);
        \filldraw[black] (-1.41,0) circle (1pt) node[anchor=west]{$r_6$};

        \draw[black, thin] (-1.17,-0.8) -- (-1.27,-0.87);
        \filldraw[black] (-1.17,-0.8) circle (1pt) node[anchor=south west]{$r_7$};

        \draw[black, thin] (-1,1) -- (-1.1,1.1);
        \draw[black, thin] (-1,1) -- (-1.1,1);
        \filldraw[black] (-1,1) circle (1pt) node[anchor=north west]{$r_8$};
    \end{tikzpicture}
    
    \caption{The first $8$ steps of the construction of a hedgehog of $S^1$.}
    \label{fig:enter-label}
\end{figure}

\begin{remark}
    From now on, we assume that for a given continuum $X$ and a sequence of continua $(X_n)$ with a sequence of points $(x_n)$ where $x_n \in X_n$, $\Hed(X,(X_n),(x_n))$ always denotes the continuum obtained by the above-mentioned construction. Moreover, we preserve the same notation for all sets and points as in the construction. 
\end{remark}

\begin{remark}
    In Construction \ref{construction:hedgehog} we use embeddings $h: X \to Q$ and $h_n: X_n \to Q_n \subseteq Q$. Without loss of generality, we will treat these embeddings as homeomorphisms from $X$ to $h(X)$, and from $X_n$ to $h_n(X_n)$, respectively. When we use $h^{-1}$ or $h_{n}^{-1}$ we understand these functions as homeomorphisms from $h(X)$ to $X$, and from $h_n(X_n)$ to $X_n$.
\end{remark}

\begin{observation}\label{obs:cutpoint}
    From Construction \ref{construction:hedgehog} it follows that $h(X) \cap h_n(X_n) = \{\bm{r}_n\}$ for each positive integer $n$. Consequently, each $\bm{r}_n$ is a cut point of $\Hed(X,(X_n),(x_n))$.
\end{observation}

\begin{theorem}\label{thm:conti}
    For each continuum $X$, for each sequence of continua $(X_n)$ and each sequence of points $(x_n)$, where $x_n \in X_n$ for all positive integers $n$, $\Hed(X,(X_n),(x_n))$ is also a continuum.
\end{theorem}
\begin{proof}
    By the construction, it follows that $\Hed(X,(X_n),(x_n))$ is a non-empty connected metric space. It remains to prove that $\Hed(X,(X_n),(x_n))$ is compact.

    Let $(y_n)$ be an arbitrary sequence in $\Hed(X,(X_n),(x_n))$ with the limit $y$. We distinguish the following cases.
    \begin{enumerate}[label=\textbf{Case \arabic*:}, leftmargin=1.25cm]
        \item Suppose there is a positive integer $m$ such that $y_n \in h_m(X_m)$ for infinitely many positive integers $n$. Then $$\lim_{n\to \infty} y_n = y \in h_m(X_m) \subseteq \Hed(X,(X_n),(x_n)),$$
        since $h_m(X_m)$ is compact.
        \item If $y_n \in h(X)$ for infinitely many positive integers $n$, then 
        $$\lim_{n\to \infty} y_n  = y \in h(X) \subseteq \Hed(X,(X_n),(x_n)),$$
        since $h(X)$ is compact.
        \item Suppose that for each n, $y_n \in h_{m(n)}(X_{m(n)})$, where $m(n_1) \neq m(n_2)$ for all but finitely many positive integers $n_1 \neq n_2$. Since $$\lim_{n\to \infty} \diam(h_{m(n)}(X_{m(n)})) = 0,$$ it follows $$\lim_{n\to \infty} h_{m(n)}(X_{m(n)}) \subseteq h(X)$$ and therefore $$\lim_{n\to \infty}  y_n = y \in h(X) \subseteq \Hed(X,(X_n),(x_n)).$$
        \item If there are infinitely many positive integers $n$ such that $y_n \in h(X)$, and infinitely many positive integers $n$ such that $y_n \notin h(X)$, it follows that the limit $\displaystyle{\lim_{n\to \infty} y_n =y\in h(X)}$ since $h_m(X_m) \cap h(X) \subseteq h(X)$, for all positive integers $m$.
    \end{enumerate}
    \qed
\end{proof}


\begin{observation}\label{obs:hedgehog_sim_star}
    Let $X$ be an arbitrary continuum, let $(X_n)$ be a sequence of Cook continua such that $X_{n+1}$ is a proper subcontinuum of $X_n$ for each positive integer $n$, $\lim\limits_{n\to\infty} \diam(X_n)=0$, and let $x_n\in\bigcap\limits_{m=1}^\infty X_m$ for each positive integer $n$. Let $\sim$ be an equivalence relation on $\Hed(X,(X_n),(x_n))$ defined by 
    $$x\sim y \ \Longleftrightarrow\ x=y \ \vee\ x,y\in h(X).$$ Then $\Hed(X,(X_n),(x_n))/_\sim$ is homeomorphic to a rigid star of Cook continua.
    
\end{observation}

\begin{lemma}\label{lemma:cook_no_arc}
If $X$ is a Cook continuum then $X$ does not contain an arc.
\end{lemma}

\begin{proof}
    Suppose, that a Cook continuum $X$ contains an arc $L\subseteq X$. Let $f:L\to [0,1]$ be a homeomorphism and let $g:[0,1]\to[0,1]$ be defined by $g(x)=1-x$. Then $f^{-1}\circ g\circ f:L\to X$ is a continuous function, which is not a constant function nor an identity -- a contradiction. \qed
\end{proof}

\begin{lemma}\label{lemma:Cook_lok_peanov}
    Let $X$ be a Cook continuum and $f:[0,1]\to X$ continuous function. Then $f$ is a constant function.
\end{lemma}

\begin{proof}

Suppose $f:[0,1]\to X$ is not a constant function.
Then there is $Y\subseteq X$ such that $f([0,1])=Y$ and $Y$ is not a singleton. By Theorem \ref{Hahn}, $Y$ is a non-degenerate Peano continuum. By \cite[Proposition 8.23]{nadler}, it follows that $Y$ contains an arc. Consequently, $X$ contains an arc, which is, by Lemma \ref{lemma:cook_no_arc}, in contradiction with $X$ being a Cook continuum.
\qed
\end{proof}
The following corollary follows directly from Lemma \ref{lemma:Cook_lok_peanov}.
\begin{corollary}\label{coro:cook_path_connected}
    The only path connected subcontinua of Cook continua are singletons.
\end{corollary}

\begin{corollary}\label{coro:Peano_Cook_konstantna}
    Let $X$ be any Cook continuum, $Y$ be any Peano continuum, and let $f:Y\to X$ be a continuous function. Then $f$ is a constant function.
\end{corollary}

\begin{proof}
    Let $g:[0,1]\to Y$ be a surjective continuous function. Note that such function exists by Theorem \ref{Hahn}. Then the function $f\circ g:[0,1]\to X$ is continuous. By Lemma \ref{lemma:Cook_lok_peanov}, $f\circ g$ is a constant function. Therefore, the function $f$ is a constant function since $g$ is surjective.
    \qed
\end{proof}

The following corollary follows directly from Lemma \ref{lemma:cook_no_arc} and the fact that a continuous image of a path connected space is a path connected space.

\begin{corollary}\label{coro:Peano_unija_Cook_konstantna}
    Let $X_1,X_2,\ldots X_n$ be Cook continua, $X = \bigcup_{i=1}^{n} X_i$, $Y$ be any Peano continuum, and let $f:Y\to X$ be a continuous function. Then $f$ is a constant function.
\end{corollary}

\section{From rigid to $0$-rigid Hedgehog}\label{sec:rig->hed}

In this section for a given non-degenerate Peano continuum X, we give an uncountable family $\mathcal{F}$ of pairwise non-homeomorphic rigid, for each positive integer $m \geq 2$, $\frac{1}{m}$- rigid, and $0$-rigid continua, respectively, such that for each $Y\in \mathcal{F}$ it holds that $Y$ contains a homeomorphic copy of $X$.
We use Cook continua and modify the construction from the previous section. 

\subsection{Rigid hedgehog}

\begin{theorem}\label{thm:peano_hedgehog_rigid}
    Let $X$ be an arbitrary Peano continuum, let $(X_n)$ be a sequence of Cook continua such that $X_{n+1}$ is a proper subcontinuum of $X_n$ for each positive integer $n$ and $\lim\limits_{n\to\infty}\diam(X_n)=0$, let $\{x_1\}=\bigcap_{n=1}^\infty X_n$ and $x_n=x_1$ for each positive integer $n$. Then $\Hed(X,(X_n),(x_n))$ is rigid.
\end{theorem}

\begin{figure}[H]
    \centering
    \begin{tikzpicture}
        \draw[black, thin] (0,0) -- (8,0);
        \filldraw[black] (4,-0.3) circle (0pt) node[anchor=north]{$h(X)$};

        \draw[black, step it] (1,0) -- (1,2);
        \filldraw[black] (1,0) circle (1pt) node[anchor=north]{$r_1$};

        \draw[black, step it] (3,0) -- (3,1.5);
        \filldraw[black] (3,0) circle (1pt) node[anchor=north]{$r_2$};

        \draw[black, step it] (2,0) -- (2,1.15);
        \filldraw[black] (2,0) circle (1pt) node[anchor=north]{$r_3$};

        \draw[black, step it] (5,0) -- (5,1);
        \filldraw[black] (5,0) circle (1pt) node[anchor=north]{$r_4$};

        \draw[black, step it] (7,0) -- (7,0.8);
        \filldraw[black] (7,0) circle (1pt) node[anchor=north]{$r_5$};

        \draw[black, step it] (4,0) -- (4,0.6);
        \filldraw[black] (4,0) circle (1pt) node[anchor=north]{$r_6$};

        \draw[black, step it] (1.6,0) -- (1.6,0.5);
        \filldraw[black] (1.6,0) circle (1pt) node[anchor=north]{$r_7$};

        \draw[black, step it] (0.4,0) -- (0.4,0.4);
        \filldraw[black] (0.4,0) circle (1pt) node[anchor=north]{$r_8$};

        \draw[black, step it] (7.8,0) -- (7.8,0.3);
        \filldraw[black] (7.8,0) circle (1pt) node[anchor=north]{$r_9$};

        \draw[black, step it] (5.3,0) -- (5.3,0.25);
        \filldraw[black] (5.3,0) circle (1pt) node[anchor=north]{$r_{10}$};

        \draw[black, step it] (1.3,0) -- (1.3,0.2);
        \filldraw[black] (1.3,0) circle (1pt) node[anchor=north]{$r_{11}$};

        \draw[black, step it] (6.6,0) -- (6.6,0.15);
        \filldraw[black] (6.6,0) circle (1pt) node[anchor=north]{$r_{12}$};
        
    \end{tikzpicture}
    
    \caption{Visualisation of the first $12$ steps of the construction of a hedgehog of an arc $X$ satisfying Theorem \ref{thm:peano_hedgehog_rigid}.}
    \label{fig:rigid}
\end{figure}
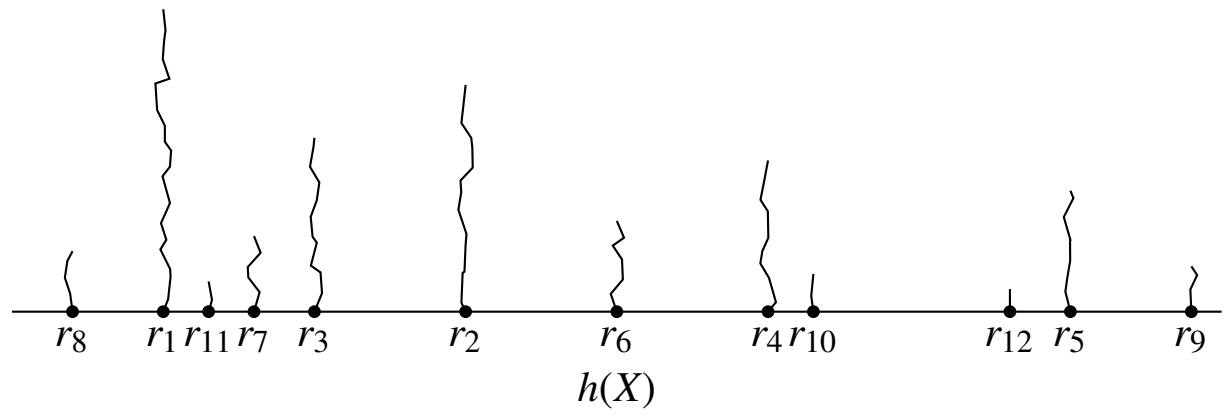

\begin{proof}
    Let $f:\Hed(X,(X_n),(x_n))\to \Hed(X,(X_n),(x_n))$ be a continuous surjective function.

    First, we prove that $f(h(X))\subseteq h(X)$. Suppose that $f(h(X))\not\subseteq h(X)$. Let $a\in h(X)$ and let $m$ be a positive integer, such that $f(a)\in h_m(X_m)\setminus h(X)$. 
    We define a function $g:\Hed(X,(X_n),(x_n))\to h_m(X_m)$ by $$g(x)=\begin{cases}
    x&  \text{;} \quad x\in h_m(X_m),\\
    h_m(x_m) &  \text{;}\quad x \notin h_m(X_m),
    \end{cases}$$
    for each $x\in\Hed(X,(X_n),(x_n))$.
    Then $g$ is a continuous surjection. Since $h(X)$ is a Peano continuum, by Corollary \ref{coro:Peano_Cook_konstantna}, $g\circ f|_{h(X)}:h(X)\to h_m(X_m)$ is a constant function that maps all $x\in h(X)$ to $f(a)$.
    Therefore, $f|_{h(X)}(x)=f(a)$ for all $x\in h(X)$.

    Let $\sim$ be an equivalence relation on $\Hed(X,(X_n),(x_n))$ defined by 
    $$x\sim y \ \Longleftrightarrow\ x=y \ \vee\ x,y\in h(X).$$
    Now, let $ \overline{f}:\Hed(X,(X_n),(x_n))/_\sim\to\Hed(X,(X_n),(x_n))/_\sim$ be a function defined by $\overline{f}([x])=[f(x)]$, which is well-defined and continuous. By Observation \ref{obs:hedgehog_sim_star} and Theorem \ref{thm:rigid_star_is_rigid}, $\Hed(X,(X_n),(x_n))/_\sim$ is a rigid continuum. 

    Since $f$ is a surjective function and $f|_{h(X)}(x)=f(a)$ for all $x\in h(X)$, $\overline{f}$ is also a surjective function. 
    Therefore, $\overline{f}=1_{\Hed(X,(X_n),(x_n))/_\sim}$ with $\overline{f}([a])=[f(a)]\neq [a]$--a contradiction.    
    Thus, $f(h(X))\subseteq h(X)$.

    Let $\sim$ and $\overline{f}$ again be the equivalence relation and the function described above. Again, by Theorem \ref{thm:rigid_star_is_rigid}, $\overline{f}=1_{\Hed(X,(X_n),(x_n))/_\sim}$.
    It follows that $f(x)=x$ for each $x\in \Hed(X,(X_n),(x_n))\setminus h(X)$. Consequently, $f(h_n(x_n))=h_n(x_n)$ for each positive integer $n$. Since $\{h_n(x_n)\ |\ n\in\NN\}=R$ is a dense subset of $h(X)$, it follows that $f(x)=x$ for each $x\in h(X)$. Hence, $f=1_{\Hed(X,(X_n),(x_n))}$, and therefore $\Hed(X,(X_n),(x_n))$ is rigid.
\qed
\end{proof}

\begin{theorem}\label{thm:rigid_uncountable}
    Let $X$ be any non-degenerate Peano continuum. There exists an uncountable family $\mathcal{F}$ of pairwise non-homeomorphic rigid continua, such that for each $Y\in \mathcal{F}$ it holds that $Y$ contains a homeomorphic copy of $X$.
\end{theorem}

\begin{proof}
    It is a well-known fact that there is a family of uncountably many pairwise non-homeomorphic Cook
continua. Denote the family of these continua by ${\cal C}$. For each $X_1 \in {\cal C}$ we can construct sequences $(x_n)$ and $(X_n)$ satisfying the properties stated in Theorem  \ref{thm:peano_hedgehog_rigid}. Consequently, $\Hed(X,(X_n),(x_n))$ is a rigid continuum. Since $X_1$ was chosen arbitrarily, the theorem follows.
    \qed
\end{proof}

\subsection{$\frac{1}{m}$-rigid hedgehog}

In the next definition, we give a special type of continua, named m-pyramids, which are obtained by Cook continua.

\begin{definition}\label{def:pyramid}
    Let $m$ be a positive integer, $m\geq 2$, and let $C$ be a Cook continuum. Let $C_1$ and $C_2$ be non-degenerate subcontinua of $C$ such that $C_1 \cap C_2 = \emptyset$. Let $c_1,d_1 \in C_1$, $c_1 \neq d_1$, and 
    $c_2,d_2 \in C_2$, $c_2 \neq d_2$. 

    The $m$-pyramid 
    $$
    \blacktriangle((C_1,c_1,d_1),(C_2,c_2,d_2),m)
    $$
of $(C_1,c_1,d_1)$ and $(C_2,c_2,d_2)$ is defined to be the quotient 
$$
\left(\textstyle \bigcup_{k=1}^{m}(C_1\times \{k\}) \cup(C_2 \times \{k\})\right)\big{/}_{\sim},
$$
where $\sim $ is the equivalence relation on $\textstyle \bigcup_{k=1}^{m}(C_1\times \{k\}) \cup(C_2 \times \{k\})$ defined by
\begin{itemize}
    \item $(d_1,1)\sim (c_1,2)\sim (d_2,1)$, $(d_1,2)\sim (c_1,3)\sim (d_2,2)$, $ \ldots $ , \\
    $(d_1,m-1)\sim (c_1,m)\sim (d_2,m-1)$, $(d_1,m)\sim (c_1,1) \sim (d_2,m)$,
    \item $(c_2,1) \sim (c_2,2) \sim \ldots \sim (c_2,m)$, and 
    \item $(y,k)\sim (y,k)$ for each $k\in\{1,2,\ldots, m\}$ and for each $y\in C_1 \cup C_2$.
\end{itemize}
\end{definition}

\begin{notation}
    We will use the following notation regarding the $m$-pyramids.
    \begin{itemize}
        \item Let $\sim$ be the equivalence relation as described in Definition \ref{def:pyramid}.
        The function $q_\sim:\bigcup_{k=1}^{m}(C_1\times \{k\}) \cup(C_2 \times \{k\})\to \blacktriangle((C_1,c_1,d_1),(C_2,c_2,d_2),m)$ is the quotient function defined by $$q_\sim((x,k))=[(x,k)]$$ for each $x\in C_1\cup C_2$ and each $k\in\{1,2,\ldots,m\}$.
        \item For each $k\in\{1,2,\ldots,m\}$, let $Y_k = q_\sim (C_1\times \{k\})$ and $W_k = q_\sim (C_2\times \{k\})$.
        \item For each $k\in\{1,2,\ldots,m\}$, let $y_k=q_\sim (c_1,k)$.
        \item $v=q_\sim (c_2,1)$.
        \item For each $k\in\{1,2,\ldots,m\}$, let $\alpha_k:C_1\to Y_k$ be a function defined by $\alpha_k(x)=q_\sim(x,k)$ for each $x\in C_1$ and let $\beta_k:C_2\to W_k$ be a function defined by $\beta_k(x)=q_\sim(x,k)$ for each $x\in C_2$. Note that $\alpha_k$ and $\beta_k$ are homeomorphisms.
    \end{itemize}
\end{notation}

\begin{figure}[H]
    \centering
    \begin{tikzpicture}
        \draw[black, thin] (-1,5) -- (6,5) node[anchor=west] {$C$};
        \draw[red, thick] (0,5) -- node[anchor=south] {$C_1$} (2,5);
        \draw[red, thick] (3,5) -- node[anchor=south] {$C_2$} (5,5);

        \filldraw[black] (0.5,5) circle (1pt) node[anchor=south] {$c_1$};
        \filldraw[black] (1.5,5) circle (1pt) node[anchor=south] {$d_1$};
        \filldraw[black] (3.5,5) circle (1pt) node[anchor=south] {$c_2$};
        \filldraw[black] (4.5,5) circle (1pt) node[anchor=south] {$d_2$};

        \draw[black, thick] (0,4) -- node[anchor=south] {$C_1\times \{1\}$} (2,4);
        \draw[black, thick] (3,4) -- node[anchor=south] {$C_2\times \{1\}$} (5,4);
        \draw[black, thick] (0,3) -- node[anchor=south] {$C_1\times \{2\}$} (2,3);
        \draw[black, thick] (3,3) -- node[anchor=south] {$C_2\times \{2\}$} (5,3);
        \filldraw[black] (1,2.7) circle (1pt);
        \filldraw[black] (1,2.4) circle (1pt);
        \filldraw[black] (1,2.1) circle (1pt);
        \filldraw[black] (4,2.7) circle (1pt);
        \filldraw[black] (4,2.4) circle (1pt);
        \filldraw[black] (4,2.1) circle (1pt);
        \draw[black, thick] (0,1.8) -- node[anchor=north] {$C_1\times \{m\}$} (2,1.8) ;
        \draw[black, thick] (3,1.8) -- node[anchor=north] {$C_2\times \{m\}$} (5,1.8) ;

        \draw [-stealth](5.5,3) -- node[anchor=south] {$q_\sim$} (6.5,3);

        \draw[black, thick] (7,1.5) -- node[anchor=east] {$W_1$} (8.5,4.5);
        \draw[black, thick] (8.5,1.2) -- node[anchor=north east] {$W_2$} (8.5,4.5);
        \draw[black, thick] (10,1.5) -- node[anchor=north east] {$W_3$} (8.5,4.5);
        
        \filldraw[black] (10,3.7) circle (1pt);
        \filldraw[black] (9.7,3.5) circle (1pt);
        \filldraw[black] (9.4,3.3) circle (1pt);

        \draw[black, thick] (11,3.5) -- node[anchor=south] {$W_m$} (8.5,4.5);

        \draw[black, thick] (7,1.5) -- node[anchor=north] {$Y_2$} (8.5,1.2);
        \draw[black, thick] (8.5,1.2) -- node[anchor=north] {$Y_3$} (10,1.5);
        \draw[black, thick] (10,1.5) -- node[anchor=north west] {$Y_4$} (10.6,2);
        \draw[black, thick] (11,2.7) -- node[anchor=north west] {$Y_m$} (11,3.5);

        \filldraw[black] (10.7,2.15) circle (1pt);
        \filldraw[black] (10.8,2.35) circle (1pt);
        \filldraw[black] (10.9,2.55) circle (1pt);

        \draw[black, densely dashed] (11,3.5) -- node[anchor=north west] {$Y_1$} (7,1.5);

        \node at (7,1.3) {$y_1$};
        \node at (8.5,1) {$y_2$};
        \node at (10,1.3) {$y_3$};
        \node at (11.2,3.5) {$y_m$};
        \node at (8.5,4.7) {$v$};

        \node at (9.1,0.5) {$\blacktriangle((C_1,c_1,d_1),(C_2,c_2,d_2),m)$};

    \end{tikzpicture}
    \caption{The $m$-pyramid as a quotient space of Cook continua.}
    \label{fig:enter-label}
\end{figure}
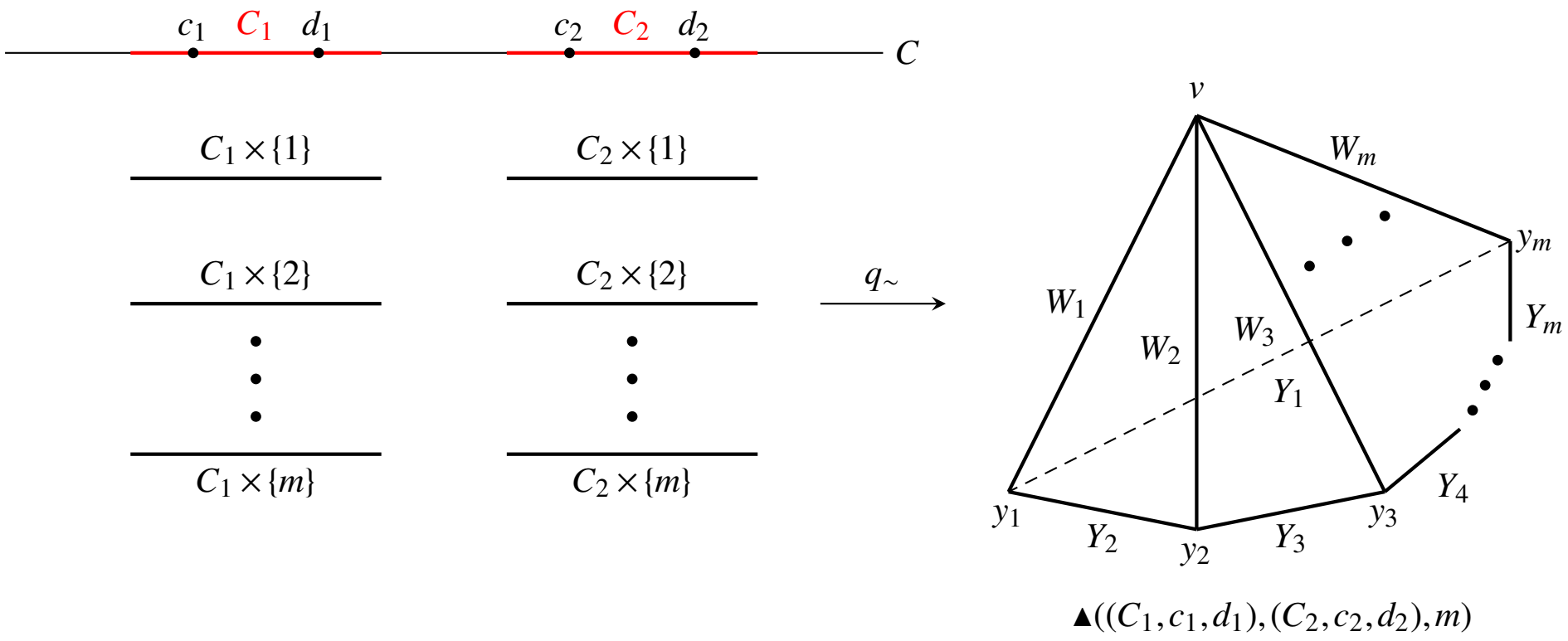

\begin{lemma}\label{lemma:m-rigid}
Let $m \geq 2$ be a positive integer and let $C$ be a Cook continuum. Let $C_1$ and $C_2$ be non-degenerate subcontinua of $C$ such that $C_1 \cap C_2 = \emptyset$. Let $c_1,d_1 \in C_1$, $c_1 \neq d_1$, and let  $c_2,d_2 \in C_2$, $c_2 \neq d_2$. Then $$\degrig(\blacktriangle((C_1,c_1,d_1),(C_2,c_2,d_2),m))=m.$$    
\end{lemma}

\begin{proof}
    Let $f:\blacktriangle((C_1,c_1,d_1),(C_2,c_2,d_2),m)\to \blacktriangle((C_1,c_1,d_1),(C_2,c_2,d_2),m)
    $ be a continuous surjective function. Let $k\in\{1,2,\ldots,m\}$ and let $Z$ be any non-degenerate subcontinuum of $W_k$. Then there exists a non-degenerate subcontinuum $B$ of the $m$-pyramid $\blacktriangle((C_1,c_1,d_1),(C_2,c_2,d_2),m)$ such that $f(B)\subseteq Z$ and $|f(B)|>1$. We prove that $B\subseteq W_1\cup W_2\cup\ldots\cup W_m$.

    Suppose, $B\cap (Y_1\cup Y_2\cup\ldots\cup Y_m\setminus\{y_1,y_2,\ldots,y_m\}) \neq \emptyset$. Without loss of generality, we assume that $B\subseteq Y_l$ for some $l\in\{1,2,\ldots,m\}$. We define a function $F:\alpha_l^{-1}(B)\to C$ by $$F(x)=\beta_k^{-1}(f(\alpha_l(x)))$$ for each $x\in \alpha_l^{-1}(B)$. The function $F$ is well-defined and continuous, and since $|f(B)|>1$, it follows that $|F(\alpha_l^{-1}(B))|>1$, and therefore, $F$ is not a constant function. Because $\alpha_l^{-1}(B)$ is a non-degenerate subcontinuum of $C_1\subseteq C$, it hold that $F(x)=x$ for each $x\in \alpha_l^{-1}(B)$. On the other hand, $F(\alpha_l^{-1}(B))\subseteq C_2$, which is a contradiction. Thus, 
    $$W_1\cup W_2\cup\ldots\cup W_m \subseteq f(W_1\cup W_2\cup\ldots\cup W_m)$$ and 
    $$(W_1\cup W_2\cup\ldots\cup W_m) \cap f((Y_1\cup Y_2\cup\ldots\cup Y_m)\setminus\{y_1,y_2,\ldots,y_m\})=\emptyset.$$

    Similarly we prove
    $$Y_1\cup Y_2\cup \ldots\cup Y_m \subseteq f(Y_1\cup Y_2\cup \ldots\cup Y_m)$$ and 
    $$(Y_1\cup Y_2\cup\ldots\cup Y_m) \cap f(W_1\cup W_2\cup\ldots\cup W_m\setminus\{y_1,y_2,\ldots,y_m\})=\emptyset.$$ 
    Consequently, $$f(W_1\cup W_2\cup\ldots\cup W_m)=W_1\cup W_2\cup\ldots\cup W_m$$ 
    and 
    $$f(Y_1\cup Y_2\cup \ldots\cup Y_m)=Y_1\cup Y_2\cup \ldots\cup Y_m.$$

    One can easily see, that $W_1\cup W_2 \cup \ldots\cup W_m$ is a star $\star_{k=1}^m(W_k,v)$ from \cite[Definition 3.3]{banicKac}, and it satisfies all the conditions from \cite[Theorem 3.7]{banicKac}. Therefore, for each $k\in\{1,2,\ldots,m\}$ there exists $l\in\{1,2,\ldots,m\}$ such that $f|_{W_k}(x)=\beta_l(\beta_k^{-1}(x))$ for each $x\in W_k$.

    The continuum $Y_1\cup Y_2\cup \ldots \cup Y_m$ is a cycle $\circ_{k=1}^m(Y_k,y_{k-1},y_k)$, where $y_0=y_m$, from \cite[Definition 3.12]{banicKac}, and it satisfies all the conditions from \cite[Theorem 3.15]{banicKac}. Therefore, for each $k\in\{1,2,\ldots,m\}$ there exists $l\in\{1,2,\ldots,m\}$ such that $f|_{Y_k}(x)=\alpha_l(\alpha_k^{-1}(x))$ for each $x\in Y_k$.
    Moreover, by \cite[Theorem 4.6]{banicKac}, the degree of rigidity of $Y_1\cup Y_2\cup \ldots \cup Y_m$ is $m$. Since $f$ is continuous, for each $k \in \{1,2,\ldots,m\}$ there is $l \in \{1,2,\ldots,m\}$ such that $f(Y_k) = Y_l$ and therefore $f(W_k) = W_l$. Consequently, the degree of rigidity of $\blacktriangle((C_1,c_1,d_1),(C_2,c_2,d_2),m)$ is also $m$.

\qed
\end{proof}

The next result follows directly from Corollary \ref{coro:Peano_unija_Cook_konstantna}.

\begin{corollary}\label{coro:piramida_peanov_konstantna}
    Let $m\geq 2$ be a positive integer and let $C$ be a Cook continuum. Let $C_1$ and $C_2$ be disjoint non-degenerate subcontinua of $C$. Let $c_1,d_1 \in C_1$, $c_1 \neq d_1$, and $c_2,d_2 \in C_2$, $c_2 \neq d_2$, and let $P= \blacktriangle((C_1,c_1,d_1),(C_2,c_2,d_2),m)$ be a $m$-pyramid. Let $X$ be any Peano continuum and let $f:X\to P$ be a continuous function. Then $f$ is a constant function.
\end{corollary}

\begin{theorem}\label{thm:1/m-rigid}
    Let $m\geq 2$ be any integer and let $C$ be a Cook continuum. Suppose $C_1, C_2, C_3$ are pairwise disjoint non-degenerate subcontinua of $C$. Let $c_1,d_1 \in C_1$, $c_1 \neq d_1$, and $c_2,d_2 \in C_2$, $c_2 \neq d_2$, and let $X_1= \blacktriangle((C_1,c_1,d_1),(C_2,c_2,d_2),m)$ be an $m$-pyramid. Let $(X_n)_{n\geq 2}$ be a sequence of Cook continua such that $X_2 = C_3$, $X_{n+1} \subsetneq X_n$, and $\lim \limits_{n\to \infty} \diam(X_n) = 0$. Denote $\{c\} = \bigcap_{n=2}^{\infty} X_n$. Let $x_1 = v$ and $x_n = c$ for each positive integer $n \geq 2$. Let $X$ be any Peano continuum. Then $\Hed( X,(X_n),(x_n)) $ is an $\frac{1}{m}$-rigid continuum and $\textup{degrig}(\Hed( X,(X_n),(x_n)))=m$.  
\end{theorem}

\begin{figure}[H]
    \centering
    \begin{tikzpicture}
        \draw[black, thin] (0,0) -- (8,0);
        \filldraw[black] (4,-0.3) circle (0pt) node[anchor=north]{$h(X)$};

        \draw[black, step it] (1,0) -- (0,2);
        \draw[black, step it] (1,0) -- (1.2,1.7);
        \draw[black, step it] (1,0) -- (2,2);
        \draw[black, step it] (0,2) -- (1.2,1.7);
        \draw[black, step it] (1.2,1.7) -- (2,2);
        \draw[black, step it] (0,2) arc (120:60:2);
        
        \filldraw[black] (1,0) circle (1pt) node[anchor=north]{$r_1$};

        \draw[black, step it] (3,0) -- (3,1.5);
        \filldraw[black] (3,0) circle (1pt) node[anchor=north]{$r_2$};

        \draw[black, step it] (2,0) -- (2,1.15);
        \filldraw[black] (2,0) circle (1pt) node[anchor=north]{$r_3$};

        \draw[black, step it] (5,0) -- (5,1);
        \filldraw[black] (5,0) circle (1pt) node[anchor=north]{$r_4$};

        \draw[black, step it] (7,0) -- (7,0.8);
        \filldraw[black] (7,0) circle (1pt) node[anchor=north]{$r_5$};

        \draw[black, step it] (4,0) -- (4,0.6);
        \filldraw[black] (4,0) circle (1pt) node[anchor=north]{$r_6$};

        \draw[black, step it] (1.6,0) -- (1.6,0.5);
        \filldraw[black] (1.6,0) circle (1pt) node[anchor=north]{$r_7$};

        \draw[black, step it] (0.4,0) -- (0.4,0.4);
        \filldraw[black] (0.4,0) circle (1pt) node[anchor=north]{$r_8$};

        \draw[black, step it] (7.8,0) -- (7.8,0.3);
        \filldraw[black] (7.8,0) circle (1pt) node[anchor=north]{$r_9$};

        \draw[black, step it] (5.3,0) -- (5.3,0.25);
        \filldraw[black] (5.3,0) circle (1pt) node[anchor=north]{$r_{10}$};

        \draw[black, step it] (1.3,0) -- (1.3,0.2);
        \filldraw[black] (1.3,0) circle (1pt) node[anchor=north]{$r_{11}$};

        \draw[black, step it] (6.6,0) -- (6.6,0.15);
        \filldraw[black] (6.6,0) circle (1pt) node[anchor=north]{$r_{12}$};
        
    \end{tikzpicture}
    
    \caption{A visualisation of the first $12$ steps of the construction of a hedgehog of an arc $X$ satisfying Theorem \ref{thm:1/m-rigid} where $m=3$.}
    \label{fig:1_m_rigid}
\end{figure}
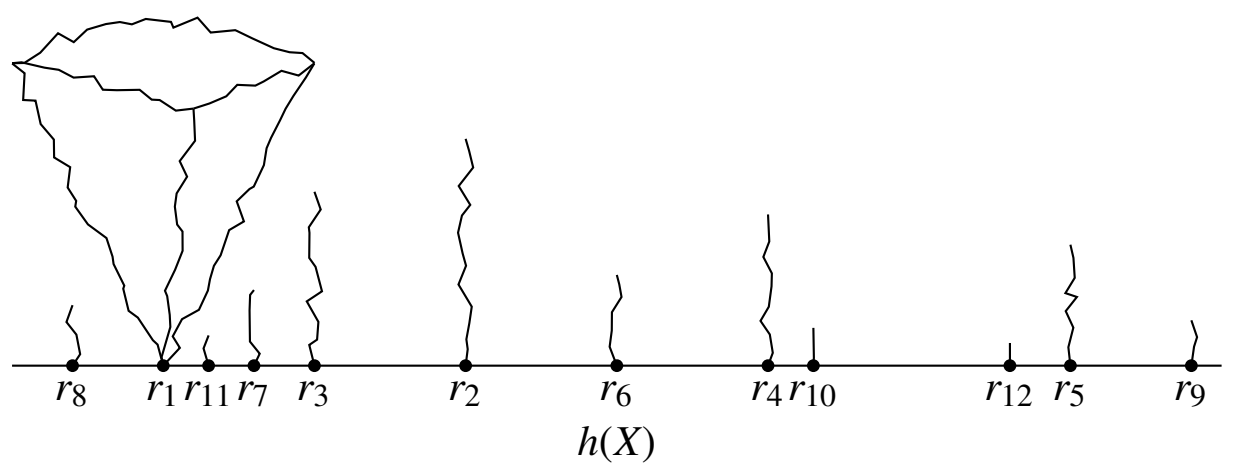

\begin{proof}
    By Theorem \ref{thm:conti}, $\Hed(X,(X_n),(x_n))$ is a continuum.
    
    Let $f: \Hed(X,(X_n),(x_n)) \to \Hed(X,(X_n),(x_n))$ be a continuous surjective function.
     For each $i \in \{ 1,2, \ldots, m\}$ we define a function $G_i:\Hed(X,(X_n),(x_n)) \to h_1(W_i)$ by 
        $$
        G_i(x)=\begin{cases}
            x & \text{;} \quad x\in h_1(W_i),\\
            (h_1\circ \beta_i\circ \beta_k^{-1}\circ h_1^{-1})(x) & \text{;} \quad x\in h_1(W_k),~ k\in\{1,2,\ldots,m\},~ k\neq i,\\
            h_1(x_1) & \text{;} \quad x\notin h_1(X_1),\\
            h_1(y_i) & \text{;} \quad x\in h_1(Y_k),~ k\in\{1,2,\ldots,m\},
        \end{cases}
        $$
        for each $x\in \Hed(X,(X_n),(x_n))$. Then $G_i$ is a continuous surjective function for each $i$.
    
    First, we prove that $f(h(X))\subseteq h(X)$. Suppose $f(h(X))\not\subseteq h(X)$. Then, there exists a point $a\in h(X)$ and there exists a positive integer $\ell$ such that $f(a)\in h_\ell(X_\ell)\setminus h(X)$. Fix such $a$ and $\ell$. The functions $g_k:\Hed(X,(X_n),(x_n))\to h_k(X_k)$, defined by 
    $$g_k(x)=\begin{cases}
        x & \text{;} \quad x\in h_k(X_k),\\
        h_k(x_k) & \text{;} \quad x\notin h_k(X_k),
    \end{cases}
    $$
    for each $x\in \Hed(X,(X_n),(x_n))$, are continuous surjections for each positive integer $k$. From Corollaries \ref{coro:Peano_Cook_konstantna} and \ref{coro:piramida_peanov_konstantna} it follows that $g_\ell\circ f|_{h(X)}$ is a constant function that maps $h(X)$ to $\{f(a)\}$. From the definition of $g_\ell$ it also holds that $f(h(X))=\{f(a)\}$.


    We distinguish the following cases.
    \begin{enumerate}[label=\textbf{Case \arabic*:}, leftmargin=1.25cm]
        \item\label{thm-proof:peano-piramid-Case1} Suppose $f(a) \in h_1(W_i)$ for some $i \in \{1,2,\ldots,m\}$.         
        For each positive integer $k\geq 2$ it holds that $G\circ f|_{h_k(X_k)}$ is a continuous function that maps $h_k(X_k)$ to $h_1(W_i)$, and therefore, it is a constant function whose image is $\{f(a)\}$. From the definition of $G$ it follows that $f(h_k(X_k))=\{f(a)\}$.
        Similarly, for each $k\in \{1,2,\ldots,m\}$ by considering functions $G_i\circ f|_{h_1(W_k)}$ and $G_i\circ f|_{h_1(Y_k)}$ we obtain that $f(h_1(W_k))=f(h_1(Y_k))=\{f(a)\}$. It follows that $f$ is a constant function--a contradiction.
        
        \item\label{proof:piramida:case2} Suppose $f(a) \in h_1(Y_i)$ for some $i \in \{1,2,\ldots,m\}$.         
        Let $k$ be a positive integer. If $f(h_k(X_k))\not\subseteq h_1(Y_1\cup Y_2\cup\ldots\cup Y_m)$, then there exists $x\in h_k(X_k)$ such that $f(x)\in h_1(W_1\cup W_2\cup\ldots\cup W_m)$ since $f(h_k(x_k))=f(a)$ and $h_1(x_1)$ is a cut point for $\Hed(X,(X_n),(x_n))$.
        Then $G_i\circ f|_{h_k(X_k)}$ is a constant function whose image is $\{h_1(y_i)\}$.
        Thus $f(h_k(X_k))\subseteq h_1(Y_1\cup Y_2\cup\ldots\cup Y_m)$.
        Consequently, $f(\Hed(X,(X_n),(x_n)))\subseteq h_1(Y_1\cup Y_2\cup \ldots \cup Y_m)$--a contradiction.

        \item Suppose $\ell \geq 2$.
        For each $i\in \{ 1,2, \ldots , m\}$, $W_i$, $Y_i$ and $X_\ell$ are disjoint subcontinua of Cook continuum $C$. Therefore, $g_\ell\circ f|_{h_1(W_i)}$ and $g_\ell\circ f|_{h_1(Y_i)}$ are constant functions for each $i\in \{ 1,2, \ldots , m\}$.     
        From $f(h_1(x_1))=\{f(a)\}$ it follows that $f(h_1(W_i))=\{f(a)\}$, and consequently $f(h_1(Y_1))=\{f(a)\}$. Hence, $f(h_1(X_1))=\{f(a)\}$.
        
        Similarly, $g_\ell \circ f|_{h_j(X_j)}$ are constant functions for each positive integer $j$. 
        From $f(h_j(x_j))=f(a)$ it also follows that $g_\ell \circ f|_{h_j(X_j)}(h_j(X_j))=\{f(a)\}$ and therefore, by definition of $g_\ell$, $f|_{h_j(X_j)}(h_j(X_j))=\{f(a)\}$ for each positive integer $j$.
        Hence, $f$ is a constant function--a contradiction.

    \end{enumerate}

    Since all three cases lead to contradictions, we conclude that $f(h(X))\subseteq h(X)$.

    Let $k\geq 2$ be any positive integer. We prove that $f(x)=x$ for each $x\in h_k(X_k)$. 
    
    Suppose $\ell$ is a positive integer with $\ell \geq 2$, $\ell\neq k$. A function $h_k^{-1}\circ g_k\circ f\circ h_{\ell}$ is a continuous function from Cook continuum $X_{\ell}$ to Cook continuum $X_k$.
    
    For $\ell<k$ it holds $X_k\subsetneqq X_{\ell}$, and therefore 
    $h_k^{-1}\circ g_k\circ f\circ h_{\ell}$ is a constant function, by Proposition \ref{prop1}. Moreover, $g_k\circ f|_{h_{\ell}(X_{\ell})}$ is a constant function, more precisely, $f(h_{\ell}(X_{\ell}))=\{h_k(x_k)\}$.

    For $\ell>k$ the function $h_k^{-1}\circ g_k\circ f\circ h_m$ is a continuous function. Since $f(h(X))\subseteq h(X)$ it holds that either $(h_k^{-1}\circ g_k\circ f\circ h_{\ell})(x)=h_k(x_k)$ for each $x\in X_{\ell}$ or $(h_k^{-1}\circ g_k\circ f\circ h_{\ell})(x)=x$ for each $x\in X_{\ell}$. 
    In either way, $f(h_{\ell}(X_{\ell}))\cap h_k(X_n\setminus X_{k+1})=\emptyset$.

    Using a similar approach as in the second case above, one can see that $f$ maps $h_1(X_1)$ to at most one point in $h_k(X_k)$. 
    As $f$ is a surjective function, something must be mapped to $h_k(X_k\setminus X_{k+1})$. The only remaining possibility is $h_k(X_k)$.
    Considering the continuous function $g_k\circ f|_{h_k(X_k)}$, we get $f(x)=x$ for each $x\in h_k(X_k)$. 
    
    Since $k$ was arbitrarily chosen, $f(x)=x$ for each $x\in \bigcup_{n=2}^\infty h_n(X_n)$, and since the set $\{h_n(x_n)~|~n \in \mathbb{N} \}$ is dense in $h(X)$ it follows that $f(x)=x$ for each $x\in  h(X)\cup \bigcup_{n=2}^\infty h_n(X_n)$.
    Moreover, $h_1(X_1)\subseteq f(h_1(X_1))$.

    Now, we consider a function $g_1\circ f|_{h_1(X_1)}$, which is continuous and surjective. Since it maps from an $m$-pyramid to itself, by Lemma \ref{lemma:m-rigid} it follows that $(g_1\circ f|_{h_1(X_1)})^m=1_{h_1(X_1)}$, and thus $f(h_1(X_1))=h_1(X_1)$.
    Consequently, $f^m=1_{\Hed(X,(X_n),(x_n))}$, and therefore $\degrig(\Hed(X,(X_n),(x_n)))\leq m$.
    It remains to prove that $\degrig(\Hed(X,(X_n),(x_n)))\geq m$.
    In order to do this, we define a function $F:\Hed(X,(X_n),(x_n))\to \Hed(X,(X_n),(x_n))$ as follows
    $$F(x)=\begin{cases}
        x & ; \quad x\in \text{Hed}(X)\setminus h_1(X_1),\\
        h_1(\alpha_{k(\Mod m) + 1}(\alpha_k^{-1}(h_1^{-1}(x)))) & ; \quad  x\in h_1(Y_k),\, k \in\{1,2,\ldots, m\},\\
        h_1(\beta_{k(\Mod m) + 1}(\beta_k^{-1}(h_1^{-1}(x))))   &;\quad x\in h_1(W_k),\, k \in\{1,2,\ldots, m\},
    \end{cases}$$
    for each $x \in \Hed(X,(X_n),(x_n))$. It follows that $F$ is a continuous surjective function and
    $F^m(x) = x $ for each $x \in h_1(X_1)$ and $F^k(x) \neq x $ for each $x \in h_1(X_1) \setminus \{x_1\}$ and each $k\in \{1,2,\ldots, m-1\}$. Therefore, the theorem follows. 
\qed
\end{proof}

\begin{corollary}
    For each Peano continuum $X$ and each positive integer $m$, there exists a $\frac{1}{m}$-rigid continuum $Y$ such that $Y$ contains a homeomorphic copy of $X$.
\end{corollary}
\begin{proof}
    The corollary follows directly from Theorem \ref{thm:1/m-rigid}.
    \qed
\end{proof}

Using similar reasoning as in the proof of Theorem \ref{thm:rigid_uncountable}, one can also prove Theorem \ref{thm:1/n_rigid_uncountable}.
 
\begin{theorem}\label{thm:1/n_rigid_uncountable}
    Let $X$ be any non-degenerate Peano continuum and let $m\geq 2$ be a positive integer. There exists an uncountable family $\mathcal{F}$ of pairwise non-homeomorphic $\frac{1}{m}$-rigid continua such that for each $Y\in \mathcal{F}$ it holds that $Y$ contains a homeomorphic copy of $X$.
\end{theorem}

\subsection{$0$-rigid hedgehog}

\begin{theorem}\label{thm:0-rigid}
    Let $C$ be a Cook continuum, and let $C_1$ and $C_2$ be pairwise disjoint non-degenerate subcontinua of $C$.
    Let $X_1$ be a simple fan of Cook continua obtained from $C_1$ as in \cite[Construction 1]{crepnjakKac} and let $(X_n)_{n\geq 2}$ be a sequence of Cook continua such that $X_2 = C_2$, $X_{n+1} \subsetneq X_n$, and $\lim \limits_{n\to \infty} \diam(X_n) = 0$. Denote $\{c\} = \bigcap_{n=2}^{\infty} X_n$. Let $x_1 = x_0$, the cut point in $X_1$, and $x_n = c$ for each $n \geq 2$. Let $X$ be any non-degenerate Peano continuum. Then the following holds.
    \begin{enumerate}
        \item $\text{Hed}( X,(X_n),(x_n)) $ is a $0$-rigid continuum.
        \item If $f: \text{Hed}( X,(X_n),(x_n)) \to \text{Hed}( X,(X_n),(x_n))$ is a continuous surjective function, then $f$ is a homeomorphism.    
    \end{enumerate}    
\end{theorem}

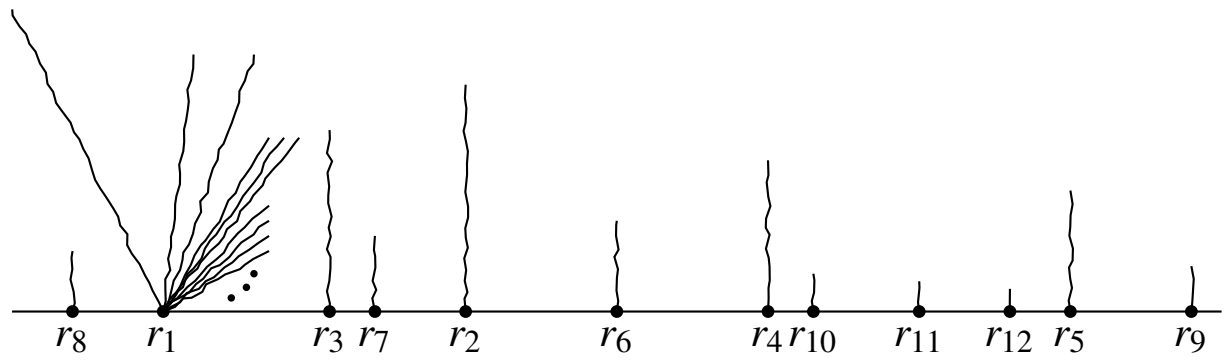
\begin{figure}[H]
    \centering
    \begin{tikzpicture}
        \draw[black, thin] (0,0) -- (8,0);
        \filldraw[black] (4,-0.3) circle (0pt) node[anchor=north]{$h(X)$};

        \draw[black, star it] (1,0) -- (0,2);
        
        \draw[black, star it] (1,0) -- (1.2,1.7);
        \draw[black, star it] (1,0) -- (1.6,1.7);
        
        \draw[black, star it] (1,0) -- (1.7,1.15);
        \draw[black, star it] (1,0) -- (1.8,1.15);
        \draw[black, star it] (1,0) -- (1.9,1.15);
        
        \draw[black, star it] (1,0) -- (1.7,0.7);
        \draw[black, star it] (1,0) -- (1.7,0.6);
        \draw[black, star it] (1,0) -- (1.7,0.5);
        \draw[black, star it] (1,0) -- (1.7,0.4);

        \filldraw[black] (1.6,0.25) circle (0.5pt);
        \filldraw[black] (1.55,0.16) circle (0.5pt);
        \filldraw[black] (1.45,0.09) circle (0.5pt);
        
        \filldraw[black] (1,0) circle (1pt) node[anchor=north]{$r_1$};

        \draw[black, star it] (3,0) -- (3,1.5);
        \filldraw[black] (3,0) circle (1pt) node[anchor=north]{$r_2$};

        \draw[black, star it] (2.1,0) -- (2.1,1.2);
        \filldraw[black] (2.1,0) circle (1pt) node[anchor=north]{$r_3$};

        \draw[black, star it] (5,0) -- (5,1);
        \filldraw[black] (5,0) circle (1pt) node[anchor=north]{$r_4$};

        \draw[black, star it] (7,0) -- (7,0.8);
        \filldraw[black] (7,0) circle (1pt) node[anchor=north]{$r_5$};

        \draw[black, star it] (4,0) -- (4,0.6);
        \filldraw[black] (4,0) circle (1pt) node[anchor=north]{$r_6$};

        \draw[black, star it] (2.4,0) -- (2.4,0.5);
        \filldraw[black] (2.4,0) circle (1pt) node[anchor=north]{$r_7$};

        \draw[black, star it] (0.4,0) -- (0.4,0.4);
        \filldraw[black] (0.4,0) circle (1pt) node[anchor=north]{$r_8$};

        \draw[black, star it] (7.8,0) -- (7.8,0.3);
        \filldraw[black] (7.8,0) circle (1pt) node[anchor=north]{$r_9$};

        \draw[black, star it] (5.3,0) -- (5.3,0.25);
        \filldraw[black] (5.3,0) circle (1pt) node[anchor=north]{$r_{10}$};

        \draw[black, star it] (6,0) -- (6,0.2);
        \filldraw[black] (6,0) circle (1pt) node[anchor=north]{$r_{11}$};

        \draw[black, star it] (6.6,0) -- (6.6,0.15);
        \filldraw[black] (6.6,0) circle (1pt) node[anchor=north]{$r_{12}$};
        
    \end{tikzpicture}
    
    \caption{A visualisation of the first $12$ steps of the construction of a hedgehog of an arc $X$ satisfying Theorem \ref{thm:0-rigid}.}
    \label{fig:1_m_rigid}
\end{figure}

\begin{proof}
   By Theorem \ref{thm:conti}, $\Hed(X)=\Hed(X,(X_n),(x_n))$ is a continuum. 

  Since $X_1$ is a simple fan of Cook continua, it holds that $X_1=\bigcup_{n=1}^\infty \bigcup_{m=1}^n K_n^m$, where $K_n^m$ are Cook continua such that 
  \begin{itemize}
      \item  $K_n^m\cap K_i^j = \{x_1\}$ if and only if $(n,m)\neq (i,j)$ for each $n,i \in \mathbb{N}$, each $m\in\{1,2,\ldots,n\}$, and for each $j\in \{1,2,\ldots, i\}$;
 \item $K_n^m$ is homeomorphic to $K_n^k$ for each $n\in \mathbb{N}$ and for each $ m,k\in \{1,2,\ldots,n\}$. Let $h_n^{m,k}:K_n^m\to K_n^k $ and $h_n^{k,m}:K_n^k\to K_n^m$ be the homeomorphisms. Note that $h_n^{m,k}\circ h_n^{k,m}$ and $h_n^{k,m}\circ h_n^{m,k}$ are the identity functions on $K_n^k$ and $K_n^m$, respectively, since all homeomorphisms between Cook continua are uniquely determined, see \cite[Lemma 2.12]{snoha}.
    \end{itemize} 
     
    
    Let $f: \text{Hed}(X) \to \text{Hed}(X)$ be a continuous surjective function. For each positive integer $m$ we define a function $g_m:\Hed(X)\to h_m(X_m)$ by
    $$
    g_m(x)=\begin{cases}
        x & \text{;} \quad x\in h_m(X_m),\\
        h_m(x_m) & \text{;} \quad x\notin h_m(X_m),
    \end{cases}
    $$
    for each $x\in \Hed(X)$. Such functions are continuous and surjective.

    First, we prove that $f(h(X))\subseteq h(X)$. Suppose that there exists $a\in h(X)$ such that $f(a)\notin h(X)$. Let $n$ be a positive integer such that $f(a)\in h_n(X_n)$. 
    We distinguish between the following two cases.

    \begin{enumerate}[label=\textbf{Case \arabic*:}, leftmargin=1.25cm]
        \item $n\neq 1$. The function $g_n\circ f|_{h(X)}$ is a continuous function which maps a Peano continuum to a Cook continuum. By Corollary \ref{coro:Peano_Cook_konstantna}, $g_n\circ f|_{h(X)}$ is a constant function, for which $(g_n\circ f|_{h(X)})(x)=f(a)$ for each $x\in h(X)$. From the definition of the function $g$, it follows that $f(x)=f(a)$ for each $x\in h(X)$.

    Let $m\geq 2$ be a positive integer. The function $g_n\circ f|_{h_m(X_m)}$ is a continuous function that maps a Cook continuum $h_m(X_m)$ to a Cook continuum $h_n(X_n)$. Since $h_m(x_m)\in h(X)$, it follows that $(g_n\circ f|_{h_m(X_m)})(h_m(x_m))=f(a)$. Therefore, $g_n\circ f|_{h_m(X_m)}$ is a constant function. Moreover, $f(x)=f(a)$ for each $x\in h_m(X_m)$.

    Since $f$ is a continuous surjection and $f(h(X)\cup \bigcup_{m=2}^\infty h_m(X_m))=\{f(a)\}$, it means that $f(h_1(X_1))=\Hed(X)$, and thus $f|_{h_1(X_1)}:h_1(X_1)\to \Hed(X)$ is a surjection. It follows that $g_1\circ f|_{h_1(X_1)}:h_1(X_1)\to h_1(X_1)$ is a continuous surjection. By \cite[Corollary 4]{crepnjakKac}, the function $g_1\circ f|_{h_1(X_1)}$ is a homeomorphism.
    Consequently, by the definition of $g_1$, $f(h_1(X_1))=h_1(X_1)$--a contradiction since $f$ is surjective.
    
        \item  $n=1$. Let $i\in \mathbb{N}$ and $j\in\{1,2,\ldots,i\}$ such that $f(a)\in h_1(K_i^j)$. We define $g_i^j:\Hed(X)\to h_1(K_i^j)$ by 
    $$
    g_i^j(x)=\begin{cases}
        x & \text{;} \quad x\in h_1(K_i^j),\\
        h_1(x_1) & \text{;} \quad x\notin h_1(K_i^j),
    \end{cases}
    $$
    for each $x\in \Hed(X)$. 
    The function $g_i^j$ is a continuous surjection. It follows that $g_i^j\circ f|_{h(X)}: h(X) \to h_1(K_i^j)$ is a continuous function which maps a Peano continuum to a Cook continuum. By Corollary \ref{coro:Peano_Cook_konstantna}, $g_i^j\circ f|_{h(X)}$ is a constant function, for which $(g_i^j\circ f|_{h(X)})(x)=f(a)$ for each $x\in h(X)$. From the definition of the function $g_i^j$, it follows that $f(x)=f(a)$ for each $x\in h(X)$.
    
    Let $m\geq 2$ be a positive integer.
The function $g_i^j\circ f|_{h_m(X_m)}$ is a continuous function that maps a Cook continuum $h_m(X_m)$ to a Cook continuum $h_1(K_i^j)$. Since $h_m(x_m)\in h(X)$, it follows that $(g_i^j\circ f|_{h_m(X_m)})(h_m(x_m))=f(a)$. Therefore, $g_i^j\circ f|_{h_m(X_m)}$ is a constant function. Moreover, $f(x)=f(a)$ for each $x\in h_m(X_m)$.

Since $f$ is a continuous surjection and $f(h(X)\cup \bigcup_{m=2}^\infty h_m(X_m))=\{f(a)\}$, it means that $f(h_1(X_1))=\Hed(X)$, and thus $f|_{h_1(X_1)}:h_1(X_1)\to \Hed(X)$ is a surjection. It follows that $g_1\circ f|_{h_1(X_1)}:h_1(X_1)\to h_1(X_1)$ is a continuous surjection. By \cite[Corollary 4]{crepnjakKac}, the function $g_1\circ f|_{h_1(X_1)}$ is a homeomorphism.
 On the other hand, there exist $a,b\in h_1(X_1)$, $a\neq b$, such that $f(a),f(b)\notin h_1(X_1)$ and therefore $(g_1\circ f|_{h_1(X_1)})(a)=(g_1\circ f|_{h_1(X_1)})(b)=h_1(x_1)$--a contradiction. 
    \end{enumerate}
    Since both cases lead to contradictions, it follows that $f(h(X))\subseteq h(X)$.

    Let $n\geq 2$ be any positive integer. We prove that $f(x)=x$ for each $x\in h_n(X_n)$. 
    
    Suppose $m$ is a positive integer with $m\geq 2$, $m\neq n$. The function $h_n^{-1}\circ g_n\circ f\circ h_m$ is a continuous function from a Cook continuum $X_m$ to a Cook continuum $X_n$. 
    
    For $m<n$, it holds that $X_n\subsetneqq X_m$, and by Proposition \ref{prop1}, $h_n^{-1}\circ g_n\circ f\circ h_m$ is a constant function. Moreover, $g_n\circ f|_{h_m(X_m)}$ is a constant function. Therefore, $f$ maps $h_m(X_m)$ to at most one point in $h_n(X_n)$, more precisely to $\{h_n(x_n)\}$.

    For $m>n$, it holds that $X_m\subsetneqq X_n$, and the function $h_n^{-1}\circ g_n\circ f\circ h_m$ is a continuous function such that either $(h_n^{-1}\circ g_n\circ f\circ h_m)(x)=x_n$ for each $x\in X_m$, since $f(h(X))\subseteq h(X)$, or $(h_n^{-1}\circ g_n\circ f\circ h_m)(x)=x$ for each $x\in X_m$. In either way, $f$ does not map any point in $h_m(X_m)$ to $h_n(X_n\setminus X_{n+1})$.

    Let $i$ be a positive integer and $j\in \{1,2,\ldots,i\}$. The function $h_n^{-1}\circ g_n\circ f|_{K_i^j}$ is a continuous function from a Cook continuum $K_i^j$ to a Cook continuum $X_n$, which are homeomorphic to distinct subcontinua of the same Cook continuum. Therefore, $h_n^{-1}\circ g_n\circ f|_{K_i^j}$ is a constant function, and consequently, $f|_{K_i^j}$ is a constant function such that $f|_{K_i^j}(K_i^j)=\{h_n(x_n)\}$.

    Since $f$ is a surjective function, something is mapped to $h_n(X_n\setminus X_{n+1})$. The only remaining possibility is $h_n(X_n)$. Considering the continuous function $ g_n\circ f|_{h_n(X_n)}$ we conclude that $f(x)=x$ for each $x\in h_n(X_n)$.

    The set $\{h_n(x_n)\ |\ n\geq 2\}$ is a dense subset of $h(X)$ for which $f(h_n(x_n))=h_n(x_n)$ for each positive integer $n\geq 2$, thus $f(x)=x$ for each $x\in h(X)$. Consequently, $f(x)=x$ for each $x\in h(X)\cup\bigcup_{n=2}^\infty h_n(X_n)$. Therefore, $h_1(X_1)\subseteq f(h_1(X_1))$.

    A continuous function $h_1^{-1}\circ g_1\circ f\circ h_1 $ maps a simple fan of Cook continua $X_1$ to itself. By \cite[Corollary 4]{crepnjakKac} the function $h_1^{-1}\circ g_1\circ f\circ h_1 $ is a homeomorphism, and thus, $g_1\circ f$ is a homeomorphism. It follows that $f|_{h_1(X_1)}$ is a homeomorphism. Therefore, $f$ is a homeomorphism.

    By \cite[Theorem 3]{crepnjakKac}, $\degrig (h_1(X_1))=\infty$ and thus $\Hed( X,(X_n),(x_n)) $ is a $0$-rigid continuum.


  
    \qed
\end{proof}

\begin{corollary}\label{cor:e_0-rigid}
    For any Peano continuum $X$, there exists a $0$-rigid continuum $Y$ such that $Y$ contains a homeomorphic copy of $X$.
\end{corollary}
\begin{proof}
    The existence of such continuum is guaranteed by Theorem \ref{thm:0-rigid}.
    \qed
\end{proof} 

Using similar reasoning as in the proof of Theorem \ref{thm:rigid_uncountable} and Corollary \ref {cor:e_0-rigid}, one can also prove Theorem \ref{thm:0_rigid_uncountable}.

\begin{theorem}\label{thm:0_rigid_uncountable}
    Let $X$ be any non-degenerate Peano continuum. There exists an uncountable family $\mathcal{F}$ of pairwise non-homeomorphic $0$-rigid continua such that for each $Y\in \mathcal{F}$ it holds that  $Y$ contains a homeomorphic copy of $X$.
\end{theorem}


\section{Application to inverse limits}\label{sec:application}

The main goal of this section is to prove that for a given continuum $X$ it holds that $\varprojlim\{\text{Hed}(X),f_n\}$ is homeomorphic to $\text{Hed}(X)$ for any sequence $(f_n)$ of continuous surjections $f_n:\text{Hed}(X)\rightarrow \text{Hed}(X)$. With these results we generalize results of Bani\v c, \v Crepnjak, and Kac, \cite{banicKac,crepnjakKac}.

\begin{theorem}\label{thm:invlim_1/m}
    Let $X$ be any non-degenerate Peano continuum and $Y$ be a $\frac{1}{m}$-rigid continuum such that $Y$ contains a homeomorphic copy of $X$, where $m$ is a positive integer. Suppose $(f_n)$ is a sequence of surjective continuous functions from $Y$ to $Y$.
    Then $\varprojlim \{Y,f_n\}$ is homeomorphic to $Y$.
\end{theorem}
\begin{proof}
    By Theorems \ref{thm:peano_hedgehog_rigid} and \ref{thm:1/m-rigid} such continuum $Y$ exists. 
    From \cite[Theorem 5.6]{banicKac} the theorem follows. 
    \qed
\end{proof}

\begin{theorem}\label{thm:invlim_0}
    Let $X$ be any non-degenerate Peano continuum. There exists a $0$-rigid continuum $Y$ such that $Y$ contains a homeomorphic copy of $X$ and for each sequence of surjective continuous functions $(f_n)$ from $Y$ to $Y$, the inverse limit $\varprojlim \{Y,f_n\}$ is homeomorphic to $Y$.
\end{theorem}
\begin{proof}
    Let $Y=\text{Hed}(X)$, obtained as in Theorem \ref{thm:0-rigid}. From the same theorem follows that $Y$ is $0$-rigid continuum which contains a homeomorphic copy of $X$ and $f_n$ is a homeomorphism for each positive integer $n$. Therefore, the theorem follows. 
    \qed
\end{proof}



The next corollary is the main result of the section.
    
\begin{corollary}\label{coro:answer} 
There exists an uncountable family $\mathcal{F}$ of pairwise non-homeomorphic continua  such that for each $X\in\mathcal{F}$ it holds
\begin{enumerate}
\item\label{enu:1} $\text{degrig}(X)=\infty$,
\item\label{enu:2} $\varprojlim\{X,f_n\}$ is homeomorphic to $X$ for any sequence $(f_n)$ of continuous surjections $f_n:X\rightarrow X$. 
\end{enumerate}
\end{corollary}
\begin{proof}  
    Let $Y_1$ and $Y_2$ be any two non-degenerate non-homeomorphic continua. 
    Let $\text{Hed}(Y_1)$ and $\text{Hed}(Y_2)$ be obtained as in Theorem \ref{thm:0-rigid}. 
    Suppose there is a homeomorphism $f:\Hed(Y_1) \to \Hed(Y_2)$. Since cut points are preserved under a homemorphism, it follows that a homeomorphic copy of $Y_1$ is mapped by $f$ to a homeomorphic copy of $Y_2$, which is a contradiction.

    By Theorems \ref{thm:0-rigid} and \ref{thm:invlim_0}, $\Hed(Y_1)$ and $\Hed(Y_2)$ satisfy \ref{enu:1} and \ref{enu:2}, and since $Y_1$ and $Y_2$ were arbitrarily chosen, the theorem follows.
    \qed

    
\end{proof}

In \cite{banicKac}, the authors posed Problem 5.7:\\
Find all non-degenerate continua $X$ such that 
\begin{enumerate}
\item $\dr(X)=\infty$,
\item $\varprojlim\{X,f_n\}$ is homeomorphic to $X$ for any sequence $(f_n)$ of continuous surjections $f_n:X\rightarrow X$. 
\end{enumerate}

This problem requires the characterization of such continua. However, the results of this paper show that for each non-degenerate continuum $X$ we can obtain the corresponding continuum $Y$ satisfying all the properties of the above-mentioned problem. Therefore, such characterization seems not to be feasible.



\end{document}